\documentclass[12pt]{amsart}
\usepackage{amssymb}
\usepackage{amsfonts}
\usepackage{amsthm}
\numberwithin{equation}{section}
\theoremstyle{plain}
\newtheorem{theorem}{Theorem}[section]
\newtheorem{proposition}[theorem]{Proposition}

\newtheorem{lemma}[theorem]{Lemma}
\theoremstyle{definition}

\newtheorem*{remark}{Remark}

\newcommand{\refE}[1]{(\ref{E:#1})}
\newcommand{\refS}[1]{Section~\ref{S:#1}}

\newcommand{\refT}[1]{Theorem~\ref{T:#1}}

\newcommand{\refP}[1]{Proposition~\ref{P:#1}}

\newcommand{\refL}[1]{Lemma~\ref{L:#1}}
\newcommand{\R}{\ensuremath{\mathbb{R}}}
\newcommand{\re}{\mathrm{\;Re\;}}
\newcommand{\im}{\mathrm{\;Im\;}}
\newcommand{\C}{\ensuremath{\mathbb{C}}}
\newcommand{\N}{\ensuremath{\mathbb{N}}}

\renewcommand{\P}{\ensuremath{\mathbb{P}}}
\newcommand{\Z}{\ensuremath{\mathbb{Z}}}
\newcommand{\K}{\ensuremath{\mathbb{K}}}
\newcommand{\F}{\ensuremath{\mathbb{F}}}

\newcommand{\scp}[2]{{\langle #1,#2\rangle}}
\newcommand{\w}{\omega}
\def\cal{\mathcal}
\def\e{e^{(m)}}
\def\T{T^{(m)}}
\def\m{^{(m)}}
\def\L{{\bf L}}
\def\d{{\rm div}}
\def\B{{\cal B}}
\def\S{{\cal S}}
\begin{document}

\vspace*{-1cm}
\hspace*{\fill} Mannheimer Manuskripte 253

{\hspace*{\fill} math/0006063}

\vspace*{2cm}

\title{Identification of Berezin-Toeplitz deformation quantization}
\author[A.V. Karabegov]{Alexander V. Karabegov}
\address[Alexander V. Karabegov]{Theory Division, Yerevan Physics
Institute, Alikhanyan bros. 2, Yerevan 375036, Armenia}
\email{karabeg@uniphi.yerphi.am}
\author[M. Schlichenmaier]{Martin Schlichenmaier}
\address[Martin Schlichenmaier]{Department of Mathematics and 
  Computer Science, University of Mannheim, D7, 27 \\
         D-68131 Mannheim \\
         Germany}
\email{schlichenmaier@math.uni-mannheim.de}
\begin{abstract}
We give a complete identification of the
deformation quantization which was obtained from
the Berezin-Toeplitz quantization on an 
arbitrary
compact K\"ahler manifold. The deformation
quantization with the opposite star-product
proves to be a
differential deformation quantization with
separation of variables whose classifying form
is explicitly calculated. Its characteristic
class (which classifies star-products up to
equivalence) is obtained. The proof is based
on the microlocal description of the Szeg\"o
kernel of a strictly pseudoconvex domain given
by Boutet de Monvel and Sj\"ostrand.
\end{abstract}
\date{April 17. 2000}
\maketitle
\section{Introduction}\label{S:intro}
In the seminal work \cite{BFFLS} Bayen, Flato, Fronsdal,
Lichnerowicz and Sternheimer drew the attention of both
physical and mathematical communities to a well posed
mathematical problem of describing and classifying up to some
natural equivalence the formal associative differential
deformations of the algebra of smooth functions on a manifold. 
The deformed associative product is traditionally denoted
$\star$ and called star-product.

If the manifold carries a Poisson structure, or a symplectic
structure (i.e. a non-degenerate Poisson structure) or even
more specific if the manifold is a K\"ahler manifold with
symplectic structure coming from the K\"ahler form one
naturally asks for a deformation of the algebra of smooth
functions in the ``direction'' of the given Poisson structure.
According to \cite{BFFLS} this deformation is treated as a
quantization of the corresponding Poisson manifold.

Due to work of De Wilde and Lecomte \cite{DeWiLe},
Fedosov \cite{Fed}, and Omori, Maeda and Yoshioka
\cite{OMY} it is known that every symplectic manifold
admits a deformation quantization in this sense.  The
deformation quantizations for a fixed symplectic
structure can be classified up to equivalence by formal
power series with coefficients in two-dimensional cohomology of
the underlying manifold, see \cite{BeCaGu},
\cite{Delstar}, \cite{Fedbook}, \cite{NeTS},
\cite{Xu}.  Kontsevich \cite{Kont} showed that every
Poisson manifold admits a deformation quantization and
that the equivalence classes of deformation
quantizations on a Poisson manifold can be parametrized
by the formal deformations of the Poisson structure. 

Despite the general existence and classification theorems it is
of importance to study deformation quantization for manifolds
with additional geometric structure and ask for deformation
quantizations respecting in a certain sense this additional
structure. Examples of this additional structure are the
structure of a complex manifold or symmetries of the manifold.

Another natural question in this context is how some naturally
defined deformation quantizations fit into the classification
of all deformation quantizations. 

In this article we will deal with
K\"ahler manifolds. Quantization of
K\"ahler manifolds via symbol
algebras was considered by Berezin in
the framework of his quantization
program developed in
\cite{Berequ},\cite{Beress}. In this
program Berezin considered symbol
algebras with the symbol product
depending on a small parameter
$\hbar$ which has a prescribed
semi-classical behavior as $\hbar\to
0$. To this end he introduced the
covariant and contravariant symbols
on K\"ahler manifolds. However, in
order to study quantization via
symbol algebras on K\"ahler manifolds
he, as well as most of his
successors, was forced to consider
K\"ahler manifolds which satisfy very
restrictive analytic conditions. 
These conditions were shown to be met
by certain classes of homogeneous
K\"ahler manifolds, e.g., $\C^n$,
generalized flag manifolds, Hermitian
symmetric domains etc. The
deformation quantization obtained
from the asymptotic expansion in
$\hbar$ as $\hbar\to 0$ of the
product of Berezin's covariant
symbols on these classes of K\"ahler
manifolds was studied in a number of
papers by Moreno, Ortega-Navarro
(\cite{Mor1}, \cite{MoOr});  Cahen,
Gutt, Rawnsley (\cite{CGR2},
\cite{CGR3}, \cite{CGR4}); 
see also \cite{Kar4}. 
This deformation quantization is
differential and respects the
separation of variables into
holomorphic and anti-holomorphic ones
in the sense that left
star-multiplication (i.e. the
multiplication with respect to the
deformed product) with local
holomorphic functions is pointwise
multiplication, and right
star-multiplication with local
anti-holomorphic functions is also
point-wise multiplication, see
\refS{sep} for the precise
definition.  It was shown in
\cite{Kar1} that such deformation
quantizations "with separation of
variables" exist for every K\"ahler
manifold.  Moreover, a complete
classification (not only up to
equivalence) of all differential
deformation quantizations with
separation of variables was given.
They are parameterized by formal
closed forms of type $(1,1)$.  The
basic results are sketched in
\refS{sep} below.  Independently a
similar existence theorem was proven
by Bordemann and Waldmann \cite{BW}
along the lines of Fedosov's
construction.  The corresponding
classifying (1,1)-form was calculated
in \cite{Kar5}. Yet another
construction was given by Reshetikhin
and Takhtajan in \cite{ResTak}. They
directly derive it from Berezin's
integral formulas which are treated
formally, i.e., with the use of the
formal method of stationary phase.
The classifying form of deformation
quantization from \cite{ResTak} 
can be easily obtained by the
methods developed in this paper. 

In \cite{Eng} Engli{\v s} obtained asymptotic
expansion of Berezin transform on a quite
general class of complex domains which do not
satisfy the conditions imposed by Berezin.

For general compact K\"ahler manifolds
$(M,\w_{-1})$ which are quantizable, i.e. admit
a quantum line bundle $L$ it was shown by
Bordemann, Meinrenken and Schlichenmaier
\cite{BMS} that the correspondence between the
Berezin-Toeplitz operators and their
contravariant symbols associated to $L^m$ has
the correct semi-classical behavior as
$m\to\infty$.  Moreover, it was shown in
\cite{Schlhab},\cite{Schlbia95}, \cite{Schldef}
that it is possible to define a deformation
quantization via this correspondence. For this
purpose one can not use the product of
contravariant symbols since in general it can
not be correctly defined.

The approach of \cite{BMS} was based on the
theory of generalized Toeplitz operators due to
Boutet de Monvel and Guillemin \cite{BGTo},
which was also used by Guillemin \cite{Guisp} in
his proof of the existence of deformation
quantizations on compact symplectic manifolds.

The deformation quantization obtained in
\cite{Schlhab},\cite{Schlbia95}, which we call
the Berezin-Toeplitz deformation quantization,
is defined in a natural way related to the
complex structure. It fulfils the condition to
be `null on constants' (i.e.  $1\star g=g\star
1=g$), it is self-adjoint (i.e.  $\overline{f\star
g}=\overline{g}\star \overline{f}$), and admits a
trace of certain type (see \cite{Schldef} for
details).  

As one of the results of this article we will
show that the Berezin-Toeplitz deformation
quantization is differential and has the
property of separation of variables, though with
the roles of holomorphic and antiholomorphic
variables swapped. To comply with the
conventions of \cite{Kar1} we consider the
opposite to the Berezin-Toeplitz deformation
quantization (i.e., the deformation quantization
with the opposite star-product) which is a
deformation quantization with separation of
variables in the usual sense. 

We will show how the Berezin-Toeplitz
deformation quantization fits into the
classification scheme of \cite{Kar1}. 
Namely, we will show that the classifying formal
(1,1)-form of its opposite deformation
quantization is 
\begin{equation}\label{E:class}
   \tilde\w=-\frac
{1}{\nu}\w_{-1}+\w_{can},
\end{equation} 
where $\nu$ is the
formal parameter, $\w_{-1}$ is the
K\"ahler form we started with and
$\w_{can}$ is the closed curvature
(1,1)-form of the canonical line
bundle of $M$ with the Hermitian
fibre metric determined by the
symplectic volume. Using \cite{Kar2}
and \refE{class} we will calculate
the classifying cohomology class
(classifying up to equivalence) of
the Berezin-Toeplitz deformation
quantization. This class was first
calculated by E. Hawkins in
\cite{Haw} by K-theoretic methods
with the use of the index theorem for
deformation quantization
(\cite{Fedbook}, \cite{NeTS}).

In deformation quantization with
separation of variables an important
role is played by the formal Berezin
transform $f\mapsto I(f)$ (see
\cite{Kar3}).  In this paper we
associate to a deformation
quantization with separation of
variables also a non-associative
"formal twisted product"
$(f,g)\mapsto Q(f,g)$.  Here the
images are always in the formal power
series over the space $C^\infty(M)$.
In the compact K\"ahler case by
considering all tensor powers $L^m$
of the line bundle $L$ and with the
help of Berezin-Rawnsley's coherent
states \cite{Raw}, it is possible to
introduce for every level $m$ the
Berezin transform $I^{(m)}$ and also
some "twisted product" $Q^{(m)}$. The
key result of this article is that
the analytic asymptotic expansions of
$I^{(m)}$, resp. of $Q^{(m)}$ define
formal objects which coincide with
$I$ and $Q$ for some deformation
quantization with separation of
variables whose classifying form $\omega$
is completely determined in terms of the form
$\tilde\omega$ (\refT{central}). To
prove this we use the integral
representation of the Szeg\"o kernel
on a strictly pseudoconvex domain
obtained by Boutet de Monvel and
Sj\"ostrand in \cite{BS} and a
theorem by Zelditch \cite{Zel} based
on \cite{BS}. We also use the method
of stationary phase and introduce its
formal counterpart which we call
"formal integral". 

Since the analytic Berezin transform $I^{(m)}$ has the
asymptotics given by the formal
Berezin transform it follows
also that the former has the expansion 
\begin{equation}
I^{(m)}=\mbox{id}+\frac 1m\Delta +O(\frac 1{m^2}),
\end{equation}
where $\Delta$ is the Laplace-Beltrami operator on $M$.

It is worth mentioning that the above formal
form $\omega$ is the formal object corresponding to the 
asymptotic expansion of the pullback of the Fubini-Study
form via Kodaira embedding of $M$ into the projective space
related to $L^m$ as $m\to+\infty$. This asymptotic
expansion was obtained by Zelditch in \cite{Zel} as a 
generalization of a theorem by Tian \cite{Tian}.

The article is organized as follows.  In \refS{sep} we
recall the basic notions of deformation quantization
and the construction of the deformation quantization
with separation of variables given by a formal
deformation of a (pseudo-)K\"ahler form.

In \refS{formal} formal integrals are introduced. Certain basic
properties, like uniqueness are shown.

In \refS{symbols} the covariant and contravariant symbols are
introduced.  Using Berezin-Toeplitz operators the
transformation $I^{(m)}$ and the twisted product $Q^{(m)}$ are
introduced. Integral
formulas for them using 2-point, resp. cyclic 3-point functions
defined via the scalar product of coherent states are given. 

\refS{bertrans} contains the key result that $I^{(m)}$ and
$Q^{(m)}$ admit a well-defined asymptotic expansion and
that the formal objects corresponding to these expansions are
given by $I$ and $Q$ respectively. 

Finally in \refS{star} the Berezin-Toeplitz star product is
identified with the help of the results obtained in
\refS{bertrans}. 

\medskip

{\bf Acknowledgements.} We would like
to thank Boris Fedosov for
interesting discussions and Mirsolav
Engli{\v s} for bringing the work of
Zelditch to our attention.  A.K.
thanks the Alexander von Humboldt
foundation and the DFG for support
and the Department of Mathematics at
the University of Mannheim for a warm
hospitality.
\newpage
\section{Deformation quantizations with separation of variables}\label{S:sep}
Given a vector space $V$, we call the elements of the space of
formal Laurent series with a finite principal part
$V[\nu^{-1},\nu]]$ formal vectors. In such a way we define
formal functions, differential forms, differential operators,
etc.  However we shall often call these formal objects just
functions, operators, and so on, omitting the word formal.

Now assume that $V$ is a Hausdorff topological vector space and
$v(m),\ m\in \R,$ is a family of vectors in $V$ which admits an
asymptotic expansion as $m\to\infty,\quad v(m)\sim \sum_{r\geq
r_0} (1/m^r)v_r$, where $r_0\in\Z$. In order to associate to
such asymptotic families the corresponding formal vectors we use
the "formalizer" $\F:v(m)\mapsto \sum_{r\geq r_0} \nu^r
v_r\in V[\nu^{-1},\nu]]$.

Let $(M,\omega_{-1})$ be a real symplectic manifold of dimension
$2n$. For any open subset $U\subset M$ denote by ${\cal F}(U) 
=C^\infty(U)[\nu^{-1},\nu]]$ the space of formal smooth
complex-valued functions on $U$.  Set ${\cal F}={\cal F}(M)$.
Denote by
$\K=\C[\nu^{-1},\nu]]$ the field of formal numbers. 

A deformation quantization on $(M,\omega_{-1})$ is an
associative $\K$-algebra structure on ${\cal F}$, with the
product $\star$ (named star-product)  given for $f=\sum \nu^j
f_j,\ g=\sum \nu^k g_k\in{\cal F}$ by the following formula:
\begin{equation}\label{E:star} f\star g=\sum_r \nu^r
\sum_{i+j+k=r} C_i(f_j,g_k).  \end{equation} In \refE{star}
$C_r,\ r=0,1,\dots$, is a sequence of bilinear mappings
$C_r:C^\infty(M)\times C^\infty(M)\to C^\infty(M)$ where
$C_0(\varphi,\psi)=\varphi\psi$ and
$C_1(\varphi,\psi)-C_1(\psi,\varphi)=i\{\varphi,\psi\}$ for
$\varphi,\psi\in C^\infty(M)$ and $\{\cdot,\cdot\}$ is the
Poisson bracket corresponding to the form $\omega_{-1}$.

Two deformation quantizations $({\cal F},\star_1)$ and $({\cal
F},\star_2)$ on $(M,\omega_{-1})$ are called equivalent if
there exists an isomorphism of algebras $B:({\cal
F},\star_1)\to({\cal F},\star_2)$ of the form $B=1+\nu
B_1+\nu^2 B_2+\dots$, where $B_k$ are linear endomorphisms of
$C^\infty(M)$. 

We shall consider only those deformation quantizations for
which the unit constant 1 is the unit in the algebra $({\cal
F},\star)$.

If all $C_r,\ r\geq 0$, are local, i.e., bidifferential
operators, then the deformation quantization is called
differential. The equivalence classes of differential
deformation quantizations on $(M,\omega_{-1})$ are bijectively
parametrized by the formal cohomology classes from
$(1/i\nu)[\omega_{-1}]+H^2(M,\C[[\nu]])$.  The formal
cohomology class parametrizing a star-product $\star$ is
called the characteristic class of this star-product and
denoted $cl(\star)$. 

A differential deformation quantization can be localized on
any open subset $U\subset M$.  The corresponding star-product
on ${\cal F}(U)$ will be denoted also $\star$.

For $f,g\in{\cal F}$ denote by $L_f,R_g$ the operators of left
and right multiplication by $f,g$ respectively in the algebra
$({\cal F},\star)$, so that $L_fg=f\star g=R_gf$. The
associativity of the star-product $\star$ is equivalent to the
fact that $L_f$ commutes with $R_g$ for all $f,g\in{\cal F}$.
If a deformation quantization is differential then 
$L_f,R_g$ are formal differential operators. 

Now let $(M,\omega_{-1})$ be pseudo-K\"ahler, i.e., a complex
manifold such that the form $\omega_{-1}$ is of type (1,1) with
respect to the complex structure. We say that a differential
deformation quantization $({\cal F},\star)$ is a deformation
quantization with separation of variables if for any open
subset $U\subset M$ and any holomorphic function $a$ and
antiholomorphic function $b$ on $U$ the operators $L_a$ and
$R_b$ are the operators of point-wise multiplication by $a$ and
$b$ respectively, i.e., $L_a=a$ and $R_b=b$. 

A formal form $\omega=(1/\nu)\omega_{-1}+\omega_0+\nu\omega_1+\dots$ is
called a formal deformation of the form $(1/\nu)\omega_{-1}$ if the forms
$\omega_r,\ r\geq 0$, are closed but not necessarily nondegenerate
(1,1)-forms on $M$.

It was shown in \cite{Kar1} that all deformation quantizations with separation
of variables on a pseudo-K\"ahler manifold $(M,\omega_{-1})$ are
bijectively parametrized by the formal deformations of the form
$(1/\nu)\omega_{-1}$. 

Recall how the star-product with separation of
variables $\star$ on $M$ corresponding to the
formal form
$\omega=(1/\nu)\omega_{-1}+\omega_0+\nu\omega_1+\dots$
is constructed.  For an arbitrary contractible
coordinate chart $U\subset M$ with holomorphic
coordinates $\{z^k\}$ let
$\Phi=(1/\nu)\Phi_{-1}+\Phi_0+\nu\Phi_1+\dots$
be a formal potential of the form $\omega$ on
$U$, i.e., $\omega=-i\partial\bar\partial\Phi$
(notice that in \cite{Kar1} - \cite{Kar5} a
potential $\Phi$ of a closed (1,1)-form
$\omega$ is defined via the formula
$\omega=i\partial\bar\partial\Phi$). 

The star-product corresponding to the form
$\omega$ is such that $L_{\partial\Phi/\partial
z^k}=\partial\Phi/\partial
z^k+\partial/\partial z^k$ and
$R_{\partial\Phi/\partial\bar
z^l}=\partial\Phi/\partial\bar
z^l+\partial/\partial\bar z^l$ on $U$.  The set
${\cal L}(U)$ of all left multiplication
operators on $U$ is completely described as the
set of all formal differential operators
commuting with the point-wise multiplication
operators by antiholomorphic coordinates
$R_{\bar z^l}=\bar z^l$ and the operators
$R_{\partial\Phi/\partial\bar
z^l}=\partial\Phi/\partial\bar
z^l+\partial/\partial\bar z^l$. One can
immediately reconstruct the star-product on $U$
from the knowledge of ${\cal L}(U)$. The local
star-products agree on the intersections of the
charts and define the global star-product
$\star$ on $M$.

One can express the characteristic class
$cl(\star)$ of the star-product with separation
of variables $\star$ parametrized by the formal
form $\omega$ in terms of this form (see
\cite{Kar2}). Unfortunately, there were wrong
signs in the formula for $cl(\star)$ in
\cite{Kar2} which should be read as follows: 
\begin{equation}\label{E:coh}
cl(\star)=(1/i)([\omega]-\varepsilon/2),
\end{equation} where $\varepsilon$ is the
canonical class of the complex manifold $M$,
i.e., the first Chern class of the canonical
holomorphic line bundle on $M$. 

Given a deformation quantization with separation
of variables $({\cal F},\star)$ on the
pseudo-K\"ahler manifold $(M,\omega_{-1})$, one
can introduce the {\it formal Berezin transform}
$I$ as the unique formal differential operator on
$M$ such that for any open subset $U\subset M$,
holomorphic function $a$ and antiholomorphic
function $b$ on $U$ the relation $I(ab)=b\star
a$ holds (see \cite{Kar3}). One can check that
$I=1+\nu\Delta+\dots$, where $\Delta$ is the
Laplace-Beltrami operator corresponding to the
pseudo-K\"ahler metric on $M$. The {\it dual}
star-product $\tilde\star$ on $M$ defined for
$f,g\in {\cal F}$ by the formula $f\tilde\star
g=I^{-1}(Ig\star
If)$ is a star-product with separation of
variables on the pseudo-K\"ahler manifold
$(M,-\omega_{-1})$. For this deformation
quantization the formal Berezin transform equals
$I^{-1}$, and thus the dual to $\tilde\star$ is
again $\star$.

Denote by $\tilde\omega=-(1/\nu)\omega_{-1}
+\tilde\omega_0+\nu\tilde\omega_1+\dots$ the
formal form parametrizing the star-product
$\tilde\star$.  The opposite to the dual
star-product, $\star'=\tilde\star^{op}$, given
by the formula $f\star' g=I^{-1}(If\star Ig)$,
also defines a deformation quantization with
separation of variables on $M$ but with the
roles of holomorphic and antiholomorphic
variables swapped. Differently said, $({\cal
F},\star')$ is a deformation quantization with
separation of variables on the pseudo-K\"ahler
manifold $(\overline{M},\omega_{-1})$ where
$\overline{M}$ is the manifold $M$ with the
opposite complex structure. The formal Berezin
transform $I$ establishes an equivalence of
deformation quantizations $({\cal F},\star)$ and
$({\cal F},\star')$. 

Introduce the following non-associative operation $Q(\cdot,\cdot)$ on
${\cal F}$. For $f,g\in{\cal F}$ set $Q(f,g)=If\star Ig=I(f\star' g)=
I(g\tilde\star f)$. We shall call it formal twisted product.
The importance of the formal twisted product will be revealed
later. 

A trace density of a deformation quantization
$({\cal F},\star)$ on a symplectic manifold $M$
is a formal volume form $\mu$ on $M$ for which
the functional $\kappa(f)=\int_M f\mu,
f\in{\cal F},$ has the trace property,
$\kappa(f\star g)=\kappa(g\star f)$ for all
$f,g\in {\cal F}$ where at least one of the
functions $f,g$ has compact support.  It was
shown in \cite{Kar3} that on a local
holomorphic chart $(U,\{z^k\})$ any formal
trace density $\mu$ can be represented in the
form $c(\nu)\exp(\Phi+\Psi)dzd\bar z$, where
$c(\nu)\in \K$ is a formal constant, $dzd\bar
z=dz^1\dots dz^nd\bar z^1\dots d\bar z^n$ is
the standard volume on $U$ and
$\Phi=(1/\nu)\Phi_{-1}+\dots,
\Psi=(1/\nu)\Psi_{-1}+\dots$ are formal
potentials of the forms $\omega,\tilde\omega$
respectively such that the relations

\begin{equation}\label{E:tr}
\partial\Phi/\partial z^k=-I(\partial\Psi/\partial z^k),\
\partial\Phi/\partial\bar z^l=-I(\partial\Psi/\partial\bar
z^l),
\mbox{ and }\Phi_{-1}+\Psi_{-1}=0
\end{equation}
hold. Vice versa, any such form is a formal trace density.

\section{Formal integrals, jets, and almost analytic functions}
\label{S:formal}

Let $\phi=(1/\nu)\phi_{-1}+\phi_0+\nu\phi_1+\dots$ and
$\mu=\mu_0+\nu\mu_1+\dots$ be, respectively, a smooth
complex-valued
formal function and a smooth formal volume form on an open set
$U\subset \R^n$. Assume that $x\in U$ is a
nondegenerate critical point of the function $\phi_{-1}$ and $\mu_0$
does not vanish at $x$. We call a $\K$-linear functional $K$ on
${\cal F}(U)$ such that
\begin{itemize}
\item[(a)] $K=K_0+\nu K_1+\dots$ is a formal distribution
supported at the
point $x$;
\item [(b)] $K_0=\delta_x$ is the Dirac distribution at the
point $x$;
\item [(c)] $K(1)=1$ (normalization condition);
\item [(d)] for any vector field $\xi$ on $U$ and $f\in{\cal
F}(U)\quad 
K\bigl(\xi f+(\xi\phi+\d_\mu\xi)f\bigr)=0$,
\end{itemize}
a (normalized) {\it formal integral at the point $x$ associated
to the pair}
$(\phi,\mu)$.

It is clear from the definition that a formal integral at a
point $x$ is independent of a particular choice of the
neighborhood $U$ and is actually associated to the germs of
$(\phi,\mu)$ at $x$.  Usually we shall consider a contractible
neighborhood $U$ such that $\mu_0$ vanishes nowhere on $U$. 

We shall prove that a formal integral at the point $x$ associated to the
pair $(\phi,\mu)$ is uniquely determined. One can also show the existence
of such a formal integral, but this fact will neither be used nor proved
in what follows. 

We call two pairs $(\phi,\mu)$ and $(\phi',\mu')$ equivalent if there
exists a formal function $u=u_0+\nu u_1+\dots$ on $U$ such that
$\phi'=\phi-u,\ \mu'=e^u\mu$.

Since the expression $\xi\phi+\d_\mu\xi$ remains invariant if
we replace the pair $(\phi,\mu)$ by an equivalent one, a formal
integral is actually associated to the equivalence class of the
pair $(\phi,\mu)$.  This means that a formal integral actually
depends on the product $e^\phi \mu$ which can be thought of as
a part of the integrand of a "formal oscillatory integral".  In
the sequel it will be shown that one can directly produce
formal integrals from the method of stationary phase.

Notice that if $K$ is a formal integral associated to a pair
$(\phi,\mu)$ it is then associated to any pair
$(\phi,c(\nu)\mu)$, where $c(\nu)$ is a nonzero formal
constant.

It is easy to show that it is enough to
check condition (d) for the coordinate vector fields
$\partial/\partial x^k$ on $U$.  Moreover, if $U$ is
contractible and such that $\mu_0$ vanishes nowhere on it, one
can choose an equivalent pair of the form $(\phi', dx)$, where
$dx=dx^1\dots dx^n$ is the standard volume form.

\begin{proposition}\label{P:uniq} 
A formal integral $K=K_0+\nu K_1+\dots$ at a
point $x$, associated to a pair
$(\phi=(1/\nu)\phi_{-1}+\phi_0+\nu\phi_1+\dots,\mu)$ is uniquely
determined. 
\end{proposition}
\begin{proof} 
We assume that $K$ is defined on a coordinate chart $(U,\{x^k\}),\
\mu=dx$, and take $f\in C^\infty(U)$. Since $\d_{dx}(\partial/\partial
x^k)=0$, the last condition of the definition of a formal
integral takes the form
\begin{equation}\label{E:simp}
K\bigl(\partial f/\partial x^k+(\partial\phi/\partial x^k)f\bigr)=0. 
\end{equation}
Equating to zero the coefficient at $\nu^r,\ r\geq 0$, of the l.h.s. of
\refE{simp} we get $K_r(\partial f/\partial x^k)+\sum_{s=0}^{r+1} 
K_s\bigl((\partial\phi_{r-s} /\partial x^k)f\bigr)=0$, which can be rewritten as a
recurrent equation 
\begin{equation}\label{E:rec}
K_{r+1}\bigl((\partial\phi_{-1}/\partial x^k)f\bigr)= \mbox{r.h.s.
depending on }
K_j,\ j\leq r.
\end{equation}

Since $x$ is a nondegenerate critical point of $\phi_{-1}$, the functions
$\partial\phi_{-1}/\partial x^k$ generate the ideal of functions
vanishing at $x$. Taking into account that $K_{r+1}(1)=0$ for $r\geq 0$
we see from \refE{rec} that $K_{r+1}$ is determined uniquely. Thus the proof 
proceeds by induction. 
\end{proof}

Let $V$ be an open subset of a complex manifold $M$ and $Z$ be a
relatively closed subset of $V$. A function $f\in C^\infty(V)$ is called
almost analytic at $Z$ if $\bar\partial f$ vanishes to infinite order
there.

Two functions $f_1,f_2\in C^\infty(V)$ are called equivalent at $Z$ if
$f_1-f_2$ vanishes to infinite order there.

Consider open subsets $U\subset \R^n$ and $\tilde U\subset \C^n$
such that $U=\tilde U\cap\R^n$, and a function $f\in C^\infty(U)$. A
function $\tilde f\in C^\infty(\tilde U)$ is called an almost analytic
extension of $f$ if it is almost analytic at $U$ and $\tilde f|_U=f$.

It is well known that every $f\in C^\infty(U)$ has an
almost analytic extension uniquely determined up to equivalence.

Fix a formal deformation
$\omega=(1/\nu)\omega_{-1}+\omega_0+\nu\omega_1+\dots$ of the
form $(1/\nu)\omega_{-1}$ on a pseudo-K\"ahler manifold
$(M,\omega_{-1})$.  Consider the corresponding star-product
with separation of variables $\star$, the formal Berezin
transform $I$ and the formal twisted product $Q$ on $M$.  We
are going to show that for any point $x\in M$ the functional
$K^I_x(f)=(If)(x)$ on ${\cal F}$ and the functional $K^Q_x$ on
${\cal F}(M\times M)$ such that $K^Q_x(f\otimes g)=Q(f,g)(x)$
can be represented as formal integrals. 

Let $U\subset M$ be a contractible coordinate chart with holomorphic
coordinates $\{z^k\}$. Given a smooth function $f(z,\bar z)$ on $U$, where $U$ is
considered as the diagonal of $\tilde U=U\times\overline{U}$, one can choose its
almost analytic extension $\tilde f(z_1, \bar z_1,z_2,\bar z_2)$ on
$\tilde U$, so that $\tilde f(z,\bar z,z,\bar z)=f(z,\bar z)$. It is a
substitute of the {\it holomorphic} function $f(z_1,\bar z_2)$ on $\tilde
U$ which in general does not exist. 

Let $\Phi=(1/\nu)\Phi_{-1}+\Phi_0+\nu\Phi_1+\dots$ be a formal potential
of the form $\omega$ on $U$ and $\tilde\Phi$ its almost analytic extension
on $\tilde U$. In particular, $\tilde\Phi(x,x)=\Phi(x)$ for $x\in U$.
Introduce an analogue of the Calabi diastatic function on
$U\times U$ by the
formula $D(x,y)=\tilde\Phi(x,y)+\tilde\Phi(y,x)-\Phi(x)-\Phi(y)$.  We
shall also use the notation
$D_k(x,y)=\tilde\Phi_k(x,y)+\tilde\Phi_k(y,x)-\Phi_k(x)-\Phi_k(y)$ so that
$D=(1/\nu)D_{-1}+D_0+\nu D_1+\dots$. 

Let $\tilde\omega$ be the formal form corresponding to the dual
star-product $\tilde\star$ of the star-product $\star$. Choose a
formal potential $\Psi$ of the form $\tilde\omega$ on $U$,
satisfying equation \refE{tr}, so that $\mu_{tr}=e^{\Phi+\Psi}dzd\bar z$
is a formal trace density of the star-product $\star$ on $U$.

\begin{theorem}\label{T:i} 
For any point $x\in U$ the functional
$K^I_x(f)=(If)(x)$ on ${\cal F}(U)$ is the formal integral at $x$
associated to the pair $(\phi^x,\mu_{tr})$, where $\phi^x(y)=D(x,y)$.
\end{theorem}

\begin{remark} In the proof of the theorem we use the notion
of jet of order $N$ of a formal function $f=\sum \nu^r f_r$ at
a given point. It is also a formal object, the formal series
of jets of order $N$ of the functions $f_r$.
\end{remark}

\begin{proof} The condition that $x$ is a nondegenerate critical
point of the function $\phi^x_{-1}(y)=D_{-1}(x,y)$ directly follows
from the fact that $\Phi_{-1}$ is a potential of the non-degenerate
(1,1)-form $\omega_{-1}$. The conditions (a-c) of the
definition of
formal integral are trivially satisfied. It remains to check the
condition (d).  Replace the pair $(\phi^x,\mu_{tr})$ by the
equivalent pair
$(\phi^x+\Phi+\Psi,dzd\bar z)=
(\tilde\Phi(x,y)+\tilde\Phi(y,x)-\Phi(x)+\Psi(y),dzd\bar z)$. Put
$x=(z_0,\bar z_0),y=(z,\bar z)$. For $\xi=\partial/\partial z^k$ the
condition (d) takes the form
\[ 
I\Bigl(\partial f/\partial
z^k+(\partial/\partial z^k)\bigl(\tilde\Phi(z_0,\bar z_0,z,\bar
z)+\tilde\Phi(z,\bar z,z_0,\bar z_0)+\Psi(z,\bar z)\bigr)f\Bigr)(z_0,\bar
z_0)=0.  
\]
We shall check it by showing that

\noindent (i) $I\bigl(\partial f/\partial z^k+(\partial\Psi/\partial
z^k)f\bigr)= -If\star(\partial\Phi/\partial z^k)$; 

\noindent (ii) $I\Bigl(\bigl(\partial\tilde\Phi(z_0,\bar z_0,z,\bar
z)/\partial z^k\bigr)f\Bigr)(z_0,\bar z_0)=0$; 

\noindent (iii) $I\Bigl(\bigl(\partial\tilde\Phi(z,\bar z,z_0,\bar
z_0)/\partial z^k\bigr)f\Bigr)(z_0,\bar z_0)=
\bigl(If\star\partial\Phi/\partial z^k\bigr)(z_0,\bar z_0)$. 

First, $I\bigl(\partial f/\partial z^k+(\partial\Psi/\partial z^k)f\bigr)=
I\bigl((\partial\Psi/\partial z^k)\tilde\star f\bigr)= If\star
I(\partial\Psi/\partial z^k)= -If\star(\partial\Phi/\partial z^k)$,
which proves (i).

The function $\psi(z,\bar z)=\tilde\Phi(z_0,\bar z_0,z,\bar z)$ is
almost antiholomorphic at the point $z=z_0$. Thus, the full jet
of the function $\partial\psi/\partial z^k$ at the point $z=z_0$ is
equal to zero, which proves (ii).

The function $\theta(z,\bar z)=\partial\tilde\Phi(z,\bar z,z_0,\bar z_0)/\partial
z^k$ is almost holomorphic at the point $z=z_0$. For a holomorphic function $a$ we
have $I(af)=I(a\tilde\star f)=If\star Ia=If\star a$.  Since $I(\theta f)(z_0,\bar
z_0)$ and $(If\star\theta)(z_0,\bar z_0)$ considered modulo $\nu^N$ depend on the
jets of finite order of the functions $\theta$ and $f$ at the point $z_0$ taken
modulo $\nu^{N'}$ for sufficiently big $N'$, we can approximate $\theta$ by a formal
holomorphic function $a$ making sure that the jets of sufficiently high order of
$\theta$ and $a$ at the point $z_0$ coincide modulo $\nu^{N'}$. Then $I(\theta
f)(z_0,\bar z_0) \equiv I(af)(z_0,\bar z_0) \equiv (If\star a)(z_0,\bar z_0)
\equiv (If\star\theta)(z_0,\bar z_0)  \pmod{\nu^N}$. Since $N$ is arbitrary,
$I(\theta f)(z_0,\bar z_0)=(If\star\theta)(z_0,\bar z_0)$ identically.  The
functions $\partial\Phi/\partial z^k$ and $\theta$ have identical holomorphic parts
of jets at the point $z_0$, i.e., all the holomorphic partial derivatives (of any
order) of these functions at the point $z_0$ coincide. Since a left
star-multiplication operator of deformation quantization with separation of
variables differentiates its argument only in holomorphic directions, we get that
$(If\star\theta)(z_0,\bar z_0)=(If\star(\partial\Phi/\partial z^k))(z_0,\bar z_0)$.
This proves (iii). 

The check for $\xi=\partial/\partial\bar z^l$ is similar, which
completes the proof of the theorem. 
\end{proof}

The following lemma and theorem can be proved by the
same methods as \refT{i}. 

\begin{lemma}\label{L:zero} 
For any vector field $\xi$ on $U$ and $x\in U\quad
I\bigl(\xi_x\phi^x)(x)=0$, where $\phi^x(y)=
D(x,y)$.
\end{lemma}
($\xi_x\phi^x$ denotes differentiation of $\phi^x$
w.r.t. the parameter $x$.)

Introduce a 3-point function $T$ on $U\times U\times U$ by the formula
$T(x,y,z)=\tilde\Phi(x,y)+\tilde\Phi(y,z)+\tilde\Phi(z,x)-\Phi(x)-
\Phi(y)-\Phi(z)$.  

\begin{theorem}\label{T:q}
For any point $x\in U$ the functional $K^Q_x$ on
${\cal F}(U\times U)$ such that $K^Q_x(f\otimes g)= Q(f,g)(x)$ is the
formal integral at the point $(x,x)\in U\times U$ associated to the pair
$(\psi^x,\mu_{tr}\otimes\mu_{tr})$, where $\psi^x(y,z)=T(x,y,z)$.
\end{theorem}

\section{Covariant and contravariant symbols}\label{S:symbols}
In the rest of the paper let $(M,\omega_{-1})$ be a compact K\"ahler
manifold. Assume that there exists
a quantum line bundle $(L,h)$ on $M$, i.e., a
holomorphic hermitian line bundle with fibre metric $h$ such
that the curvature of the canonical connection on $L$ coincides with
the K\"ahler form $\omega_{-1}$.

Let $m$ be a non-negative integer.  The metric $h$ induces the
fibre metric $h^m$ on the tensor power $L^m=L^{\otimes m}$. Denote by
$L^2(L^m)$ the Hilbert space of square-integrable sections of $L^m$ with
respect to the norm $\Vert s \Vert^2=\int h^m(s) \Omega$, where
$\Omega=(1/n!)(\omega_{-1})^n$ is the symplectic volume form on $M$.  The
Bergman projector $B_m$ is the orthogonal projector in $L^2(L^m)$
onto the space $H_m=\Gamma_{hol}(L^m)$ of holomorphic sections of $L^m$. 

Denote by $k$ the metric on the dual line bundle $\tau:L^*\to
M$ induced by $h$. It is a well known fact that $D=\{\alpha\in
L^*|k(\alpha) < 1\}$ is a strictly pseudoconvex domain in
$L^*$. Its boundary $X=\{\alpha\in L^*|k(\alpha)=1\}$ is a
$S^1$-principal bundle. 

The sections of $L^m$ are identified with the $m$-homogeneous functions
on $L^*$ by means of the mapping $\gamma_m: s\mapsto \psi_s$, where
$\psi_s(\alpha)=\scp{\alpha^{\otimes m}}{s(x)}$ for $\alpha\in L^*_x$.
Here $\scp{\cdot}{\cdot}$ denotes the bilinear pairing between $(L^*)^m$
and $L^m$. 

There exists a unique $S^1$-invariant volume form
$\tilde\Omega$ on $X$ such that for every 
$f\in C^\infty(M)$ the equality
$\int_X (\tau^* f) \tilde\Omega=\int_M f\Omega$ 
holds.

The mapping $\gamma_m$ maps $L^2(L^m)$ isometrically onto the weight
subspace of $L^2(X,\tilde\Omega)$ of weight $m$ with respect to
the $S^1$-action.  The Hardy space ${\cal H}\subset
L^2(X,\tilde\Omega)$ of square integrable traces of holomorphic
functions on $L^*$ splits up into weight spaces, ${\cal
H}=\oplus_{m=0}^\infty{\cal H}_m$, where ${\cal H}_m=\gamma_m(H_m)$. 

Denote by $S$ and $\hat B_m$ the Szeg\"o and Bergman orthogonal
projections in $L^2(X,\tilde\Omega)$ onto ${\cal H}$ and ${\cal H}_m$
respectively. Thus $S=\sum_{m=0}^\infty \hat B_m$.
The Bergman projection $\hat B_m$ has a smooth integral
kernel $\B_m=\B_m(\alpha,\beta)$ on $X\times X$.

For each $\alpha\in L^*-0$ ('$-0$' means the zero section removed)
one can define a coherent state $\e_\alpha$
as the unique holomorphic section of $L^m$ such that for each
$s\in H_m\
\scp{s}{\e_\alpha}=\psi_s(\alpha)$ where
$\scp{\cdot}{\cdot}$ is the hermitian scalar product on $L^2(L^m)$
antilinear in the second argument. 

Since the line bundle $L$ is positive it is known that there exists
a constant $m_0$ such that for $m > m_0\ \dim H_m > 0$ and all
$\e_\alpha,\alpha\in L^*-0$, are nonzero vectors. From now on we
assume that $m > m_0$ unless otherwise specified.

The coherent state $\e_\alpha$ is antiholomorphic in $\alpha$ and
for a nonzero $c\in\C\ \e_{c\alpha}=\bar c^m\e_\alpha$.  Notice that
in \cite{CGR1} coherent states are parametrized by the points of
$L-0$.

For $s\in L^2(L^m)\ \scp{s}{\e_\alpha}= \scp{s}{B_m\e_\alpha}=
\scp{B_ms}{\e_\alpha}=\psi_{B_ms}(\alpha)$. 
 The mapping $\gamma_m$ intertwines the Bergman projectors $B_m$ and
$\hat B_m$, for $s\in L^2(L^m)\ \psi_{B_ms}=\hat B_m\psi_s$. Thus, on
the one hand, $\scp{s}{\e_\alpha}=\hat B_m\psi_s(\alpha)=\int_X
\B_m(\alpha,\beta)\psi_s(\beta)\tilde\Omega(\beta)$.  On the other
hand, $\scp{s}{\e_\alpha}=\scp{\psi_s}{\psi_{\e_\alpha}}=\int_X
\psi_s(\beta) \overline{\psi_{\e_\alpha}(\beta)}\tilde\Omega(\beta)$.
Taking into account that
$\scp{\e_\beta}{\e_\alpha}=\psi_{\e_\beta}(\alpha)=
\overline{\psi_{\e_\alpha}(\beta)}$ we finally get that
$\scp{\e_\beta}{\e_\alpha}=\psi_{\e_\beta}(\alpha)=\B_m(\alpha,\beta)$. 
In particular, one can extend the kernel $\B_m(\alpha,\beta)$ from
$X\times X$ to a holomorphic function on
$(L^*-0)\times(\overline{L^*-0})$ such that for nonzero $c,d\in\C$
\begin{equation}\label{E:homo} 
\B_m(c\alpha,d\beta)=(c\bar
d)^m\B_m(\alpha,\beta). 
\end{equation}

For $\alpha,\beta\in L^*-0$ the following inequality holds.

\begin{equation}\label{E:ineq}
\left| \B_m(\alpha,\beta)\right|=\left|\scp{\e_\alpha}{\e_\beta}\right| \leq
\|\e_\alpha\|\|\e_\beta\|=(\B_m(\alpha,\alpha)\B_m(\beta,\beta))^\frac{1}{2}.
\end{equation}

The covariant symbol of an operator $A$ in the space $H_m$ is the 
function
$\sigma(A)$ on $M$ such that
\[
\sigma(A)(x)=\frac{\scp{A\e_\alpha}{\e_\alpha}}
{\scp{\e_\alpha}{\e_\alpha}}
\]
for any $\alpha\in L^*_x-0$. 

Denote by $M_f$ the multiplication operator by a
function $f\in C^\infty(M)$ on sections of $L^m$. Define
the Berezin-Toeplitz operator $\T_f=B_mM_fB_m$ in $H_m$.
If an operator
in $H_m$ is represented in the form $\T_f$ for some
function $f\in C^\infty(M)$ then the function $f$ is
called its contravariant symbol. 

With these symbols we associate two important operations on
$C^\infty(M)$, the Berezin transform $I\m$ and a
non-associative binary operation $Q\m$ which we call twisted
product, as follows. For $f,g\in C^\infty(M)\ I\m
f=\sigma(\T_f),\ Q\m (f,g)=\sigma(\T_f\T_g)$. 

We are going to show in \refS{bertrans} that
both $I\m$ and $Q\m$ have asymptotic
expansions in $1/m$ as $m\to +\infty$, such that if the
asymptotic parameter $1/m$ in these expansions is replaced by
the formal parameter $\nu$ then we get the formal Berezin
transform $I$ and the formal twisted product $Q$ corresponding
to some deformation quantization with separation of variables
on $(M,\omega_{-1})$ which can be completely identified. We
shall mainly be interested in the opposite to its dual
deformation quantization. The goal of this paper is to show
that it coincides with the Berezin-Toeplitz deformation
quantization obtained in \cite{Schlbia95},\cite{Schldef}. 

In order to obtain the asymptotic expansions of $I\m$ and $Q\m$ we
need their integral representations.  To calculate them it is
convenient to work on $X$ rather than on $M$.  We shall use the fact
that for $f\in C^\infty(M),\ s\in \Gamma(L^m),
\psi_{M_fs}=(\tau^*f)\cdot\psi_s$. For $x\in M$ denote by $X_x$ the
fibre of the bundle $X$ over $x$, $X_x=\tau^{-1}(x)\cap X$. 
For $x,y,z\in M$ choose $\alpha\in X_x,\ \beta\in X_y,\ \gamma\in
X_z$ and set 
\begin{multline}\label{E:uvw}
u_m(x)=\B_m(\alpha,\alpha),\qquad
v_m(x,y)=\B_m(\alpha,\beta)\B_m(\beta,\alpha),\\
w_m(x,y,z)=\B_m(\alpha,\beta)\B_m(\beta,\gamma)\B_m(\gamma,\alpha).
\end{multline} 
It follows from \refE{homo} that $u_m(x),v_m(x,y),w_m(x,y,z)$ do not
depend on the choice of $\alpha,\beta,\gamma$ and thus relations
\refE{uvw} correctly define functions $u_m,v_m,w_m$.  
The function $w_m$ is the so called cyclic 3-point function
studied in \cite{BerSchlcse}. Notice that
$u_m(x)=\B_m(\alpha,\alpha)=\|\e_\alpha\|^2 > 0,\
v_m(x,y)=\B_m(\alpha,\beta)\B_m(\beta,\alpha)=|\B_m(\alpha,\beta)|^2\geq
0$ and 
\begin{equation}\label{E:func}
|w_m(x,y,z)|^2=v_m(x,y)v_m(y,z)v_m(z,x).
\end{equation}

It follows from \refE{ineq} that 
\begin{equation}\label{E:est}
   v_m(x,y)\leq u_m(x) u_m(y).
\end{equation}
 
For $\alpha\in X_x$ we have
\begin{multline}\label{E:ii}
\bigl(I\m
f\bigr)(x)=\sigma\bigl(\T_f\bigr)(x)=\frac{\scp{\T_f\e_\alpha}{\e_\alpha}}{
\scp{\e_\alpha}{\e_\alpha}}=
\frac{\scp{B_mM_fB_m\e_\alpha}{\e_\alpha}}{\B_m(\alpha,\alpha)}=\\
\frac{\scp{M_f\e_\alpha}{\e_\alpha}}{\B_m(\alpha,\alpha)}=
\frac{\scp{(\tau^*f)\psi_{\e_\alpha}}{\psi_{\e_\alpha}}}{\B_m(\alpha,\alpha)}= 
\frac{1}{\B_m(\alpha,\alpha)}\int_X (\tau^*f) 
\psi_{\e_\alpha}(\beta)\overline{\psi_{\e_\alpha}(\beta)}\tilde\Omega(\beta)=\\ 
\frac{1}{\B_m(\alpha,\alpha)}\int_X\B_m(\alpha,\beta) 
\B_m(\beta,\alpha)(\tau^*f)(\beta)\tilde\Omega(\beta)
=\frac{1}{u_m(x)}\int_M v_m(x,y) f(y)\Omega(y).
\end{multline}
  
Similarly we obtain that 
\begin{multline}\label{E:qq}
Q\m (f,g)(x)=\\
\frac{1}{\B_m(\alpha,\alpha)}
\int_{X\times X} \B_m(\alpha,\beta)
\B_m(\beta,\gamma)\B_m(\gamma,\alpha)(\tau^*f)(\beta)(\tau^*g)(\gamma)
\tilde\Omega(\beta)\tilde\Omega(\gamma)=\\
\frac{1}{u_m(x)}\int_{M\times M} w_m(x,y,z)
f(y)g(z)\Omega(y)\Omega(z).
\end{multline}

\section{Asymptotic expansion of the Berezin transform}\label{S:bertrans}

In \cite{BS} a microlocal description of the
integral kernel $\S$ of the Szeg\"o projection
$S$ was given. The results in \cite{BS} were
obtained for a strictly pseudoconvex domain with
a smooth boundary in $\C^{n+1}$. However,
according to the concluding remarks in
\cite{BS}, these results are still valid for the
domain $D$ in $L^*$ (see also \cite{BMS},
\cite{Zel}). 

It was proved in \cite{BS} that the Szeg\"o kernel $\S$ is a
generalized function on $X\times X$ singular on the diagonal of
$X\times X$ and smooth outside the diagonal. The Szeg\"o kernel
$\S$ can be expressed via the Bergman kernels $\B_m$ as
follows, $\S=\sum_{m\geq 0} \B_m$, where the sum should be
understood as a sum of generalized functions.

For $(\alpha,\beta)\in X\times X$ and $\theta\in\R$ set
$r_\theta(\alpha,\beta)=(e^{i\theta}\alpha,\beta)$.  Since
each ${\cal H}_m$ is a weight space of the $S^1$-action in
the Hardy space ${\cal H}$, one can recover $\B_m$ from the
Szeg\"o kernel, \begin{equation}\label{E:weak}
\B_m=\frac{1}{2\pi}\int_0^{2\pi} e^{-im\theta} r^*_\theta\S
d\theta.  \end{equation} This equality should be understood
in the weak sense. 

Let $E_1,E_2$ be closed disjoint subsets of $M$. Set
$F_i=\tau^{-1}(E_i)\cap X,\ i=1,2$. Thus $F_1,F_2$ are closed
disjoint subsets of $X$ or, equivalently, $F_1\times F_2$
is a closed subset of $X\times X$ which does not intersect
the diagonal. For $\S$ and $\B_m$ considered as smooth
functions outside the diagonal
of $X\times X$ equality \refE{weak} holds in
the ordinary sense, from whence it follows immediately that
\begin{equation}\label{E:vv} 
\sup_{F_1\times
F_2}|\B_m|=O\biggl(\frac{1}{m^N}\biggr)  
\end{equation}
for any $N\in\N$. 

Now let $E$ be a closed subset of $M$ and $x\in M\setminus E$. 
Then \refE{vv} implies that 
\begin{equation}\label{E:dec}
\sup_{y\in E} v_m(x,y)=O\biggl(\frac{1}{m^N}\biggr) 
\end{equation}
for any $N\in\N$.

In \cite{Zel} Zelditch proved that the function $u_m$ on $M$ expands
in the asymptotic series $u_m \sim m^n \sum_{r\geq 0} (1/m^r)b_r$ as
$m\to +\infty$, where $b_0=1$ ($n=(1/2)\dim_\R M$).  More
precisely, he proved that for any $k,N\in\N$
\begin{equation}\label{E:zel}
\left|u_m-\sum_{r=0}^{N-1}m^{n-r}b_r\right|_{C^k}=O(m^{n-N}).
\end{equation}

Therefore 
\begin{equation}\label{E:in}
\sup_M \frac{1}{u_m}=O\biggl(\frac{1}{m^n}\biggr).
\end{equation}

Using \refE{ii},\refE{dec} and
\refE{in} it is easy to prove the following proposition.

\begin{proposition}\label{P:loci} 
Let $f\in C^\infty(M)$ be a function vanishing in a neighborhood
of a point $x\in M$. Then $|(I\m f)(x)|= O(1/m^N)$ for any
$N\in\N$, i.e., $(I\m f)(x)$ is rapidly decreasing as $m\to
+\infty$. 
\end{proposition}

Thus for arbitrary $f\in C^\infty(M)$ and $x\in M$ the asymptotics of
$(I\m f)(x)$ as $m\to +\infty$ depends only on the germ of the
function $f$ at the point $x$. 

Let $E$ be a closed subset of $M$. Fix a
point $x\in M\setminus E$. The function
$w_m(x,y,z)$ with $y\in E$ can be
estimated using \refE{func} and \refE{est}
as follows.  
\begin{equation}\label{E:estim}
|w_m(x,y,z)|^2\leq v_m(x,y)u_m(x)u_m(y)\bigl(u_m(z)\bigr)^2.
\end{equation}
Using \refE{dec}, \refE{zel} and \refE{estim} we obtain that 
for any $N\in\N$
\begin{equation}\label{E:w}
\sup_{y\in
E,z\in M}\left|w_m(x,y,z)\right|=O\biggl(\frac{1}{m^N}\biggr).
\end{equation}
Similarly,
\begin{equation}\label{E:w'}
\sup_{y\in
M,z\in E}\left|w_m(x,y,z)\right|=O\biggl(\frac{1}{m^N}\biggr)
\end{equation}
for any $N\in\N$.

Using \refE{qq}, \refE{in}, \refE{w} and \refE{w'} 
one can readily prove the following proposition.

\begin{proposition}\label{P:locq} 
For $x\in M$ and arbitrary functions $f,g\in C^\infty(M)$ such
that $f$ or $g$ vanishes in a neighborhood of $x \quad 
Q\m (f,g)(x)$ is rapidly decreasing as $m\to +\infty$. 
\end{proposition}

This statement can be reformulated as follows.  For arbitrary
$f,g\in C^\infty(M)$ and $x\in M$ the asymptotics of $Q\m
(f,g)(x)$ as $m\to +\infty$ depends only on the germs of the
functions $f,g$ at the point $x$.

We are going to show how formal integrals can be obtained from the
method of stationary phase.

Let $\phi$ be a smooth function on an open subset $U\subset
M$ such that (i)  $\re\phi\leq 0$; (ii) there is only one
critical point $x_c\in U$ of the function $\phi$, which is
moreover a nondegenerate critical point; (iii) $\phi(x_c)=0$.

Consider a classical symbol $\rho(x,m)\in S^0(U\times\R)$ 
(see \cite{Hoer} for definition and notation) which has an
asymptotic expansion $\rho\sim\sum_{r\geq 0} (1/m^r)\rho_r(x)$ such that
$\rho_0(x_c)\neq 0$, and a smooth nonvanishing volume form $dx$ on 
$U$. Set $\mu(m)=\rho(m,x)dx$.
 
We can apply the method of stationary phase with a complex phase
function (see \cite{Hoer} and \cite{MS}) to the integral
\begin{equation}\label{E:osc}
S_m(f)=\int_U e^{m\phi} f \mu(m),
\end{equation}
where $f\in C^\infty_0(U)$. Notice that the phase function in
\refE{osc} is $(1/i)\phi$ so that the condition 
$\im\bigl((1/i)\phi\bigr)\geq 0$ is satisfied.

Taking into account that $\dim_\R M=2n$ and $\phi(x_c)=0$ we
obtain that $S_m(f)$ expands to an asymptotic series
$S_m(f)\sim \sum_{r=0}^\infty (1/m^{n+r}) \tilde K_r(f)$ as
$m\to +\infty$.  Here $\tilde K_r,\ r\geq 0,$ are
distributions supported at $x_c$ and $\tilde
K_0=c_n\delta_{x_c}$, where $c_n$ is a nonzero constant. 
Thus $\F(S_m(f))=\nu^n \tilde K(f)$, where $\F$ is the
"formalizer" introduced in \refS{sep} and $\tilde K$ is the
functional defined by the formula $\tilde K=\sum_{r\geq
0}\nu^r\tilde K_r$.  Consider the normalized functional
$K(f)=\tilde K(f)/\tilde K(1)$, so that $K(1)=1$.  Then
$\F(S_m(f))=c(\nu)K(f)$, where $c(\nu)=\nu^n c_n+\dots$ is a
formal constant. 

\begin{proposition}\label{P:stat}
For $f\in C_0^\infty(U)\ S_m(f)$ given by
\refE{osc} expands in an asymptotic series in
$1/m$ as $m\to+\infty$. 
$\F\bigl(S_m(f)\bigr)=c(\nu)K(f)$, where $K$ is
the formal integral at the point $x_c$
associated to the pair $((1/\nu)\phi,\F(\mu))$
and $c(\nu)$ is a nonzero formal constant. 
\end{proposition}

\begin{proof} 
Conditions (a-c) of the definition of formal integral are
satisfied. It remains to check condition (d). Let $\xi$ be a
vector field
on $U$. Denote by $\L_\xi$ the corresponding Lie derivative. We have
$0=\int_U\L_\xi(e^{m\phi}f\mu(m))= \int_Ue^{m\phi}\bigl(\xi
f+(m\xi\phi+\d_\mu \xi)f\bigr)\mu(m)$. Applying $\F$ we
obtain that
$0=\F\Bigl(\int_Ue^{m\phi}\bigl(\xi f+(m\xi \phi+\d_\mu
\xi)f\bigr)\mu(m)\Bigr)=c(\nu)K\Bigl(\xi f+\bigl(\xi((1/\nu)
\phi)+\d_{\F(\mu)}
\xi\bigr)f\Bigr)$, which concludes the proof.
\end{proof}

Our next goal is to get an asymptotic expansion
of the Bergman kernel $\B_m$ in a neighborhood
of the diagonal of $X\times X$ as $m\to+\infty$.
An asymptotic expansion of $\B_m$ on the
diagonal of $X\times X$ was obtained in
\cite{Zel} (see \refE{zel}). As in \cite{Zel},
we use the integral representation of the
Szeg\"o kernel $\S$ given by the following
theorem. We denote $n=\dim_\C M$.

\begin{theorem}\label{T:bs} {\bf(L. Boutet de Monvel and 
J. Sj\"ostrand, \cite{BS}, Theorem 1.5. and ${\cal x}$ 2.c)}
Let $\S (\alpha,\beta)$ be the Szeg\"o kernel of the boundary
$X$ of the bounded strictly pseudoconvex domain $D$ in the complex
manifold $L^*$. There exists a classical symbol $a\in
S^n(X\times X\times \R^+)$ which has an asymptotic expansion 
\[
a(\alpha,\beta,t)\sim \sum_{k=0}^\infty t^{n-k}
a_k(\alpha,\beta)
\]
so that
\begin{equation}\label{E:szeg}
\S(\alpha,\beta)=\int_0^\infty e^{it\varphi(\alpha,\beta)}
a(\alpha,\beta,t)dt,
\end{equation}
where the phase $\varphi(\alpha,\beta)\in C^\infty(L^*\times
L^*)$ is
determined by the following properties:
\begin{itemize}
\item $\varphi(\alpha,\alpha)=(1/i)\bigl(k(\alpha)-1\bigr)$;
\item $\bar\partial_\alpha\varphi$ and $\partial_\beta\varphi$
vanish to infinite order along the diagonal;
\item
$\varphi(\alpha,\beta)=\overline{-\varphi(\beta,\alpha)}$.
\end{itemize}
\end{theorem}

The phase function $\varphi$ is thus almost
analytic at the diagonal of $L^*\times
\overline{L^*}$. It is determined up to
equivalence at the diagonal. 

Fix an arbitrary point $x_0\in M$.  Let $s$ be a
local holomorphic frame of $L^\ast$ over a
contractible open neighborhood $U\subset M$ of
the point $x_0$ with local holomorphic
coordinates $\{z^k\}$. Then
$\alpha(x)=s(x)/\sqrt{k(s(x))}$ is a smooth
section of $X$ over $U$. Set $\Phi_{-1}(x)=\log
k\bigl(s(x)\bigr)$, so that
\begin{equation}\label{E:alph}
\alpha(x)=e^{(-1/2)\Phi_{-1}(x)}s(x).
\end{equation}
It follows from the fact
that $L$ is a quantum line bundle (i.e., that
$\omega_{-1}$ is the curvature form of the
Hermitian holomorphic line bundle $L$) that
$\Phi_{-1}$ is a potential of the form
$\omega_{-1}$ on $U$. 

Let $\tilde\Phi_{-1}(x,y)\in C^\infty(U\times \overline{U})$ be an almost
analytic extension of the potential $\Phi_{-1}$ from the
diagonal of $U\times \overline{U}$. Denote
$D_{-1}(x,y):=\tilde\Phi_{-1}(x,y)+\tilde\Phi_{-1}(y,x)
-\Phi_{-1}(x)-\Phi_{-1}(y)$. 
Since $\tilde\Phi_{-1}(x,x)=\Phi_{-1}(x)$, we have $D_{-1}(x,x)=0$.
In local coordinates
\begin{equation}\label{E:dandq}
D_{-1}(x,y)=-Q_{x_0}(x-y)+O(|x-y|^3), 
\end{equation}
where
\[
Q_{x_0}(z)=\sum
\frac{\partial^2\Phi_{-1}}{\partial z^k\partial
\bar z^l}(x_0)z^k\bar z^l
\] 
is a positive definite quadratic form (since $\omega_{-1}$ is 
a K\"ahler form).

The following statement is an immediate consequence of
\refE{dandq}.

\begin{lemma}\label{L:aae} 
There exists a neighborhood $U'\subset U$ of the point $x_0$
such that for any two different points $x,y\in U'$ one has
$\re D_{-1}(x,y) < 0$.
\end{lemma}

Taking, if necessary,
$(1/2)\bigl(\tilde\Phi_{-1}(x,y)+\overline{\tilde\Phi_{-1}(y,x)}\bigr)$
instead
of $\tilde\Phi_{-1}(x,y)$ choose $\tilde\Phi_{-1}$ such that
$\tilde\Phi_{-1}(y,x)=\overline{\tilde\Phi_{-1}(x,y)}$. Replace $U$ by a
smaller
neighborhood (retaining for it the notation $U$) such that $\re
D_{-1}(x,y) < 0$ for any different $x,y$ from this neighborhood.

For a point $\alpha$ in the restriction $L^*|_U$
of the line bundle $L^*$ to $U$ 
represented in the form $\alpha=v s(x)$ with
$v\in\C,x\in U$ one has $k(\alpha)=|v|^2     
k\bigl(s(x)\bigr)$. 

One can choose the
phase function $\varphi(\alpha,\beta)$ in
\refE{szeg} of the form

\begin{equation}\label{E:phas}
\varphi(\alpha,\beta)=(1/i)\bigl(v\bar
we^{\tilde\Phi_{-1}(x,y)}-1\bigr),
\end{equation}
where $\alpha=vs(x),\beta=ws(y)\in L^*|_U$.

Denote $\chi(x,y):=\tilde\Phi_{-1}(x,y)-
(1/2)\Phi_{-1}(x)-(1/2)\Phi_{-1}(y)$.
Notice that $\chi(x,x)=0$.

The following theorem is a slight generalization
of Theorem 1 from \cite{Zel}. 

\begin{theorem}\label{T:as} There exists an
asymptotic expansion of the Bergman kernel
$\B_m(\alpha(x),\alpha(y))$ on $U\times U$ as
$m\to+\infty$, of the form
\begin{equation}\label{E:aseb}
\B_m(\alpha(x),\alpha(y))\sim m^n
e^{m\chi(x,y)}\sum_{r\geq 0}(1/m^r)\tilde
b_r(x,y)
\end{equation}
such that (i) for any compact
$E\subset U\times U$ and $N\in \N$
\begin{equation}\label{E:estex}
\sup_{(x,y)\in E}
\left|\B_m(\alpha(x),\alpha(y))-m^n
e^{m\chi(x,y)}\sum_{r=0}^{N-1}(1/m^r)\tilde
b_r(x,y)\right|=O(m^{n-N});
\end{equation}

(ii) $\tilde b_r(x,y)$ is an almost analytic
extension of $b_r(x)$ from the diagonal of
$U\times U$, where $b_r,\ r\geq 0,$ are given by
\refE{zel}; in particular, $\tilde b_0(x,x)=1$.
\end{theorem}

\begin{proof} Using integral representations
\refE{weak} and \refE{szeg} one gets for $x,y\in
U$
\begin{equation}\label{E:zel1}
\B_m(\alpha(x),\alpha(y))=\frac{1}{2\pi}
\int_0^{2\pi} \int_0^\infty e^{-im\theta}
e^{it\varphi(r_\theta\alpha(x),\alpha(y))}
a(r_\theta\alpha(x),\alpha(y),t) d\theta dt.
\end{equation}
Changing variables $t\mapsto mt$ in \refE{zel1}
gives
\begin{equation}\label{E:zel2}
\B_m(\alpha(x),\alpha(y))=\frac{m}{2\pi}
\int_0^{2\pi} \int_0^\infty 
e^{im(t\varphi(r_\theta\alpha(x),
\alpha(y))-\theta)}
a(r_\theta\alpha(x),\alpha(y),mt) d\theta dt.
\end{equation}

In order to apply the method of stationary phase to the integral in
\refE{zel2} the following preparations should be made.

Using \refE{phas} and \refE{alph}
express the phase function of the integral in
\refE{zel2} as follows:
\begin{equation}\label{E:bigph}
Z(t,\theta;x,y):= t\varphi(r_\theta\alpha(x),
\alpha(y))-\theta=
(t/i)\bigl(e^{i\theta}
e^{\chi(x,y)}-1\bigr)-\theta.
\end{equation}

In order to find the critical points of the phase $Z$
(with respect to the variables $(t,\theta)$; the variables $(x,y)$
are parameters)
consider first the equation
\begin{equation}\label{E:partt}
\partial_t Z(t,\theta;x,y)=
(1/i)\bigl(e^{i\theta} e^{\chi(x,y)}-1\bigr)=0.
\end{equation}
It follows from $\tilde\Phi_{-1}(y,x)=\overline{\tilde\Phi_{-1}(x,y)}$ 
that $\re \chi(x,y)=
(1/2)D_{-1}(x,y)$. Since $D_{-1}(x,y) < 0$ for $x\neq y$
one has $|e^{\chi(x,y)}|=
e^{\re \chi(x,y)} < 1$ for $x\neq y$ whence it follows
that \refE{partt} holds only if $x=y$ and thus $Z$ has
critical points only if $x=y$. Since $\chi(x,x)=0$
one gets that
$\partial_t Z(t,\theta;x,x)=(1/i)(e^{i\theta}-1)$ and
$\partial_\theta Z(t,\theta;x,x)=te^{i\theta}-1$.
As in the proof of Theorem 1 from \cite{Zel}, one shows
that for each $x\in U$ the only critical point of the phase function
$Z(t,\theta;x,x)$ is $(t=1,\theta=0)$. It does not depend on
$x$ and, moreover, is nondegenerate.

One has $\im Z(t,\theta;x,y)= \im
\bigl((t/i)(e^{i\theta}e^{\chi(x,y)}-1)-\theta\bigr)=
t\bigl(1-\re
(e^{i\theta}e^{\chi(x,y)})\bigr)\geq 0$
since $|e^{\chi(x,y)}|\leq 1$. 

Finally, a simple calculation shows that the germs of the functions
$Z(t,\theta;x,y)$ and $(1/i)\chi(x,y)$ at the point
$(t=1,\theta=0,x=x_0,y=x_0)$ are equal modulo the ideal generated by
$\partial_t Z$ and $\partial_\theta Z$.

Applying now the method of stationary
phase to the integral in \refE{zel2} one obtains
the expansion \refE{aseb} satisfying \refE{estex}. 

It follows from \refE{zel} and \refE{uvw} that 
$\tilde b_r(x,x)=b_r(x)$ and $\tilde b_0(x,x)=b_0(x)=1$.
It remains to show that all $\tilde b_r, r\geq 0,$ are almost analytic
along the diagonal of $U\times \overline{U}$. One has 
\[
\B_m(\alpha(x),\alpha(y))=
e^{(-m/2)(\Phi_{-1}(x)+\Phi_{-1}(y))}\B_m(s(x),s(y)).
\]
The function $\B_m(s(x),s(y))$ is holomorphic on $U\times\overline{U}$.
Let $\xi$ and $\eta$ be arbitrary holomorphic and antiholomorphic vector
fields on $U$, respectively. Then $\xi_y\B_m(s(x),s(y))=0$ and
$\eta_x\B_m(s(x),s(y))=0$ (the subscripts $x,y$ show in which variable
the vector field acts). Thus
\[
\bigl(\eta_x+\frac{m}{2} \eta_x
\Phi_{-1}(x)\bigr)\B_m(\alpha(x),\alpha(y))=
e^{(-m/2)\Phi_{-1}(x)} \eta_x
e^{(m/2)\Phi_{-1}(x)}\B_m(\alpha(x),\alpha(y))=0.
\]
Analogously, $\bigl(\xi_y+(m/2) \xi_y
\Phi_{-1}(x)\bigr)\B_m(\alpha(x),\alpha(y))=0.$
Let $A_N$ be a product of $N$ derivations on $U\times U$. Then, using
integral representation \refE{zel2}, expand $0=A_N \bigl(\eta_x+(m/2)
\eta_x \Phi_{-1}(x)\bigr)\B_m(\alpha(x),\alpha(y))$ to the asymptotic
series
\begin{equation}\label{E:capan}
A_N \bigl(\eta_x+\frac{m}{2} \eta_x \Phi_{-1}(x)\bigr)\bigl(m^n
e^{m\chi(x,y)}\sum_{r\geq
0} (1/m^r)\tilde b_r(x,y)\bigr)=
e^{m\chi(x,y)}\sum_{r\geq r_0} (1/m^r) c_r(x,y)
\end{equation}
for some $c_r\in C^\infty(U\times U)$ and $r_0\in\Z$, and with the norm
estimate of the partial sums in the r.h.s. term in \refE{capan} analogous
to \refE{estex}.
Since $\chi(x,x)=0$ one gets that all $c_r(x,x)=0$. From this fact one
can prove by induction over $N$ that $\eta_x\tilde b_r$ vanishes to
infinite order at the diagonal of $U\times U$. Similarly, $\xi_y\tilde
b_r$ vanishes to infinite order at the diagonal. Thus
$\tilde b_r$ is almost analytic along the diagonal.
\end{proof}

Choose a symbol $b(x,y,m)\in S^0\bigl((U\times U)\times \R)$
such that it has the asymptotic expansion $b\sim
\sum_{r=0}^\infty (1/m^r)\tilde b_r$. Then
$\B_m(\alpha(x),\alpha(y))$ is asymptotically equivalent to $m^n
e^{m\chi(x,y)} b(x,y,m)$ on $U\times U$. One has
$\chi(x,y)+\chi(y,x)=\tilde\Phi_{-1}(x,y)+\tilde\Phi_{-1}(y,x) 
-\Phi_{-1}(x)-\Phi_{-1}(y)=D_{-1}(x,y)$ and $\chi(x,y)+\chi(y,z)+\chi(z,x)=
\tilde\Phi_{-1}(x,y)+\tilde\Phi_{-1}(y,z)+\tilde\Phi_{-1}(z,x)
-\Phi_{-1}(x)-\Phi_{-1}(y)-\Phi_{-1}(z)=T_{-1}(x,y,z)$ (the last equality
is the definition of $T_{-1}$). Thus the functions
\begin{multline*}
v_m(x,y)=\B_m(\alpha(x),\alpha(y))\B_m(\alpha(y),\alpha(x)) 
\mbox{ and } \\
w_m(x,y,z)=\B_m(\alpha(x),\alpha(y))\B_m(\alpha(y),\alpha(z))
\B_m(\alpha(z),\alpha(x)) 
\end{multline*}
are asymptotically equivalent to 
\[
m^{2n} e^{mD_{-1}(x,y)}b(x,y,m) b(y,x,m) \mbox { and }
m^{3n} e^{mT_{-1}(x,y,z)}b(x,y,m)b(y,z,m)b(z,x,m)
\]
respectively. It is easy to show that for the
functions $\phi^x_{-1}(y)=D_{-1}(x,y)$ and
$\psi^x_{-1}(y,z)=T_{-1}(x,y,z)$ the points
$y=x$ and $(y,z)=(x,x)$ respectively are
nondegenerate critical ones.

Since $\tilde b_0(x,x)=1$ one can take a smaller
contractible neighborhood $V\Subset U$ of $x_0$
such that $\tilde b_0(x,y)$ does not vanish on
the closure of $V\times V$. 
One can choose $V$ such that for any 
$x\in V$ the only critical
points of the functions $\phi^x_{-1}(y)$ on $V$ and 
$\psi^x_{-1}(y,z)$ on $V\times V$ are $y=x$ and $(y,z)=(x,x)$
respectively.

The identity
$T_{-1}(x,y,z)=(1/2)(D_{-1}(x,y)+D_{-1}(y,z)+D_{-1}(z,x))$
implies that\\ 
$\re T_{-1}(x,y,z)\leq 0$ for $x,y,z\in V$. 

The symbol $b(x,y,m)$ does not vanish on $V\times V$ 
for sufficiently big values of $m$. It follows from
\refE{zel} that $1/u_m(x)$ and $(m^n b(x,x,m))^{-1}$ are 
asymptotically equivalent for $x\in V$. 
Denote
\begin{equation}\label{E:meas}
\mu_x(m)=\frac{b(x,y,m)b(y,x,m)}{b(x,x,m)}\Omega(y),\
\tilde\mu_x(m)=\frac{b(x,y,m)b(y,z,m)b(z,x,m)}{b(x,x,m)}\Omega(y)\Omega(z). 
\end{equation}

Taking into account \refE{ii} we get for $f,g\in
C^\infty_0(V)$ and $x\in V$ the following asymptotic equivalences,
\begin{equation}\label{E:integr} 
\bigl(I\m f\bigr)(x)\sim m^n\int_V e^{m\phi^x_{-1}} f \mu_x(m)\mbox{ and }
Q\m(f,g)(x)\sim m^{2n}\int_{V\times V} e^{m\psi^x_{-1}} (f\otimes g)\tilde\mu_x(m).
\end{equation}
(In \refE{integr} $\bigl(f\otimes g\bigr)(y,z)=f(y)g(z)$.)

Applying \refP{stat} to the first integral in \refE{integr} we
obtain that $\F\bigl(\bigl(I\m
f\bigr)(x)\bigr)=c(\nu,x)L_x^I(f)$, where the functional
$L_x^I$ on ${\cal F}(V)$ is the formal integral at the point
$x$ associated to the pair $((1/\nu)\phi^x_{-1},\F(\mu_x))$ and
$c(\nu,x)$ is a formal function. It is easy to show that
$c(\nu,x)$ is smooth.

Similarly we obtain from \refE{integr} that
$\F\bigl(Q\m(f,g)(x)\bigr)=d(\nu,x)L_x^Q(f\otimes g)$ where the
functional $L_x^Q$ on ${\cal F}(V\times V)$ is the formal
integral at the point $(x,x)$ associated to the pair
$((1/\nu)\psi^x_{-1},\F(\tilde\mu_x))$ and $d(\nu,x)$ is a
smooth formal function.

Since the unit constant $1$ is a contravariant symbol of the
unit operator ${\bf 1}$, $\T_1={\bf 1}$, and $\sigma({\bf 1})=
1$, we have $I\m 1=1, Q\m (1,1)=1,$ and thus $\F(I\m 1)=1$ and
$\F\bigl(Q\m(1,1)\bigr)=1$. Taking the functions $f,g$ in
\refE{integr} to be equal to $1$ in a neighborhood of $x$ and
applying \refP{loci} and \refP{locq} we get that $c(\nu,x)=1$
and $d(\nu,x)=1$.

Since $b_0(x,y)$ does not vanish on $V\times V$ we can find a
formal function $\tilde s(x,y)$ on $V\times V$ such that
$\F(b(x,y,m))=e^{\tilde s(x,y)}$.  Set $s(x)=\tilde s(x,x)$. 
In these notations 
\begin{multline}\label{E:fmux}
\F(\mu_x)=\exp (\tilde s(x,y)+\tilde
s(y,x)-s(x))\Omega(y)\mbox{ and }\\
\F(\tilde\mu_x)=\exp (\tilde s(x,y)+\tilde
s(y,z)+\tilde s(z,x)-s(x))\Omega(y)\Omega(z). 
\end{multline}
It follows from \refT{as} that $\tilde
s$ is an almost analytic extension of the function $s$ from
the diagonal of $V\times V$.  According to \refE{zel},
$\F(u_m)=(1/\nu^n)e^s$.

Denote $\tilde \Phi=(1/\nu)\tilde\Phi_{-1}+\tilde s,\
\Phi=(1/\nu)\Phi_{-1}+s,\
D(x,y)=\tilde\Phi(x,y)+\tilde\Phi(y,x)-\Phi(x)-\Phi(y)
=(1/\nu)D_{-1}(x,y)+(\tilde s(x,y)+\tilde
s(y,x)-s(x)-s(y)),\
T(x,y,z)=\tilde\Phi(x,y)+\tilde\Phi(y,z)+\tilde\Phi(z,x) 
-\Phi(x)-\Phi(y)-\Phi(z)$. The pair
$((1/\nu)\phi^x_{-1},\F(\mu_x))=
((1/\nu)\phi^x_{-1},\exp(\tilde s(x,y)+\tilde
s(y,x)-s(x))\Omega(y))$ is then equivalent to the pair
$(\phi^x,e^s\Omega)$, where $\phi^x(y)=D(x,y)$.
Similarly, the pair
$((1/\nu)\psi^x_{-1},\F(\tilde\mu_x))$ is equivalent to
the pair $(\psi^x,e^s\Omega\otimes e^s\Omega)$, where
$\psi^x(y,z)=T(x,y,z)$.

Thus we arrive at the following proposition.

\begin{proposition}\label{P:bertr} 
For $f,g\in C_0^\infty(V),\ x\in V,\quad \bigl(I\m
f\bigr)(x)$ and $Q\m(f,g)(x)$ expand in asymptotic series
in $1/m$ as $m\to +\infty$.  $\F\bigl(\bigl(I\m
f\bigr)(x)\bigr)=L^I_x(f)$ and
$\F\bigl(Q\m(f,g)(x)\bigr)=L^Q_x(f\otimes g)$, where the
functional $L^I_x$ on ${\cal F}(V)$ is the formal
integral at the point $x$ associated to the pair
$(\phi^x,e^s\Omega)$ and the functional $L^Q_x$ on ${\cal
F}(V\times V)$ is the formal integral at the point
$(x,x)$ associated to the pair $(\psi^x,e^s\Omega\otimes
e^s\Omega)$. 
\end{proposition}

Now let $\star$ denote the star-product with
separation of variables on $(V,\omega_{-1})$
corresponding to the formal deformation
$\omega=-i\partial\bar\partial\Phi$ of the form
$(1/\nu)\omega_{-1}$, so that $\Phi$ is a formal
potential of $\omega$. Let $I$ be the
corresponding formal Berezin transform,
$\tilde\omega$ the formal form parametrizing the
dual star-product $\tilde\star$ and $\Psi$ the
solution of \refE{tr} so that
$\mu_{tr}=e^{\Phi+\Psi}dzd\bar z$ is a formal
trace density for the star-product $\star$. 

Choose a classical symbol $\rho(x,m)\in
S^0(V\times\R)$ which has an asymptotic
expansion $\rho\sim\sum_{r\geq 0} (1/m^r)\rho_r$
such that \begin{equation}\label{E:dens}
\F(\rho)e^s\Omega=\mu_{tr}. \end{equation}
Clearly, \refE{dens} determines $\F(\rho)$
uniquely. 

For $f\in C_0^\infty(V)$ and $x\in V$ consider the following
integral
\begin{equation}\label{E:formi}
(P_m f)(x)=m^n\int_V e^{m\phi^x_{-1}}f\rho\mu_x,
\end{equation}
where $\phi^x_{-1}(y)=D_{-1}(x,y)$ and $\mu_x$ is given by
\refE{meas}.

\begin{proposition}\label{P:formber}
For $f\in C_0^\infty(V)$ and $x\in V$ $(P_m f)(x)$ 
has an asymptotic expansion in $1/m$ as $m\to +\infty$.
$\F\bigl((P_m f)(x)\bigr)=c(\nu)(If)(x)$, where $c(\nu)$ is a
nonzero formal
constant.
\end{proposition}
\begin{proof}
It was already shown that the phase function $(1/i)\phi^x_{-1}$ of
integral \refE{formi} satisfies the conditions required in the
method of stationary phase. Thus \refP{stat} can be applied
to \refE{formi}. We get that 
$\F\bigl((P_m f)(x)\bigr)=c(\nu,x)K_x(f)$,
where $K_x$ is a formal integral at the point $x$ associated
to the pair
$((1/\nu)\phi^x_{-1},\F(\rho\mu_x))$ and $c(\nu,x)$ is a
nonvanishing formal function on $V$.
It follows from
\refE{fmux} and \refE{dens} that
$\F(\rho\mu_x)=\F(\rho)\F(\mu_x)=\F(\rho)\exp (\tilde
s(x,y)+\tilde
s(x,y)-s(x))\Omega(y)=\exp(\tilde s(x,y)+\tilde
s(x,y)-s(x)-s(y))\mu_{tr}=
\exp(D(x,y)-(1/\nu)D_{-1}(x,y))\mu_{tr}
=\exp(\phi^x-(1/\nu)\phi^x_{-1})\mu_{tr}$, where 
$\phi^x(y)=D(x,y)$.
The pair $((1/\nu)\phi^x_{-1},\F(\rho\mu_x))$ is thus equivalent
to the pair $(\phi^x,\mu_{tr})$. 
Applying \refT{i} we get that 
\begin{equation}\label{E:fff}
\F\bigl((P_m f)(x)\bigr)=c(\nu,x)\bigl(If\bigr)(x). 
\end{equation}
It remains to show that $c(\nu,x)$ is actually a formal
constant. Let $x_1$ be an arbitrary point of $V$. Choose
a function $\epsilon\in C_0^\infty(V)$ such that
$\epsilon=1$ in a neighborhood $W\subset V$ of $x_1$.
Let $\xi$ be a vector field on $V$. Then,
using \refE{fmux}, we obtain
\begin{multline}\label{E:psi}
\frac{1}{\nu}\xi_x\phi^x_{-1}(y)+
\F\biggl(\frac{\xi_x\mu_x}{\mu_x}(y)\biggr)=\\
\frac{1}{\nu}\xi_x D_{-1}(x,y)+\xi_x(\tilde
s(x,y)+\tilde s(x,y)-s(x))=\xi_x D(x,y)=\xi_x\phi^x.
\end{multline}
On the one hand,
taking into account \refE{psi} we get for $x\in W$ that
\begin{multline}\label{E:xi}
\F\bigl((\xi P_m\epsilon)(x)\bigr)=
\F\biggl(m^n\xi\int_V e^{m\phi^x_{-1}}\epsilon\rho\mu_x\biggr)=
\F\biggl(m^n\int_V e^{m\phi^x_{-1}}\Bigl(
m\xi_x\phi_x+\frac{\xi_x\mu_x}{\mu_x}\Bigr)
\epsilon\rho\mu_x\biggr)=\\
c(\nu,x)I\biggl(\frac{1}{\nu}\xi_x\phi^x_{-1}(y)+
\F\biggl(\frac{\xi_x\mu_x}{\mu_x}(y)\biggr)\biggr)=
c(\nu,x)I(\xi_x\phi^x)=0.
\end{multline}
The last equality in \refE{xi} follows from \refL{zero}.
On the other hand, for $x\in W$ we have from \refE{fff}
that $\F\bigl((P_m \epsilon)(x)\bigr)=c(\nu,x)$, from whence 
$\F\bigl((\xi P_m\epsilon)(x)\bigr)=
\xi\F\bigl((P_m\epsilon)(x)\bigr)=\xi c(\nu,x)$.
Thus we get from \refE{xi} that $\xi c(\nu,x)=0$ 
on $W$ for an arbitrary vector field $\xi$, from which the
Proposition follows.
\end{proof}

It follows from \refE{integr} and \refE{formi} that for
$f\in C_0^\infty(V)$ $\bigl(I\m(f\rho)\bigr)(x)$ is
asymptotically equivalent to $(P_mf)(x)$. Passing to
formal asymptotic series we get from \refP{bertr} and
\refP{formber}
that $c(\nu)(If)(x)=\F\bigl((P_m f)(x)\bigr)=
\F\bigl(\bigl(I\m(f\rho)\bigr)(x)\bigr)=
L^I_x(f\F(\rho))$, where $L^I_x$
is the formal integral at the point $x$ associated to the
pair $(\phi^x,e^s\Omega)$. Thus
\begin{equation}\label{E:fident}
c(\nu)(If)(x)=L^I_x(f\F(\rho)).
\end{equation}
The formal function $\F(\rho)$ is invertible
(see \refE{dens}). Setting $f=1/\F(\rho)$ in
\refE{fident} we get
$c(\nu)\bigl(I(1/\F(\rho))\bigr)(x)=L^I_x(1)=1$
for all $x\in V$. Since the formal Berezin
transform is invertible and $I(1)=1$, we finally
obtain that
\begin{equation}\label{E:final}
\F(\rho)=c(\nu).
\end{equation}

Now \refE{dens} can be rewritten as follows,
\begin{equation}\label{E:rel}
c(\nu)e^s\Omega=d\mu_{tr}=e^{\Phi+\Psi}dzd\bar z.
\end{equation}

In local holomorphic coordinates the
symplectic volume $\Omega$ can be expressed as
follows, $\Omega=e^\theta dzd\bar z$.  The
closed (1,1)-form
$\omega_{can}=-i\partial\bar\partial\theta$
does not depend on the choice of local
holomorphic coordinates and is defined
globally on $M$. The form $\omega_{can}$ is
the curvature form of the canonical
connection of the canonical holomorphic line
bundle on $M$ equipped with the Hermitian
fibre metric determined by the volume form
$\Omega$. Its de Rham class
$\varepsilon=[\omega_{can}]$ is the first
Chern class of the canonical holomorphic line
bundle on $M$ and thus depends only on the
complex structure on $M$.  The class
$\varepsilon$ is called the canonical class of
the complex manifold $M$.

One can see from \refE{rel} that 
$c(\nu)=c_0+\nu c_1+\dots$, where $c_0\neq 0$. 
Thus there exists a formal constant $d(\nu)$ such that
$e^{d(\nu)}=c(\nu)$ and $d(\nu)+s+\theta=\Phi+\Psi$.
Therefore the formal potential $\Psi$ of the form
$\tilde\omega$ is expressed explicitly,
$\Psi=d(\nu)-(1/\nu)\Phi_{-1}+\theta$, from whence it
follows that 
\begin{equation}\label{E:omega}
\tilde\omega=-(1/\nu)\omega_{-1}+\omega_{can}.
\end{equation} 

Formula \refE{omega} defines $\tilde\omega$ globally on
$M$. Thus the corresponding star-product $\tilde\star$
and therefore its dual star-product $\omega$ are also
globally defined. 

\refT{i}, \refT{q}, \refP{loci}, \refP{locq}
\refP{bertr},
formulas \refE{fident},
\refE{final} and \refE{rel} imply the following theorem,
which is the
central technical result of the paper. 

\begin{theorem}\label{T:central}
For any $f,g\in C^\infty(M)$ and $x\in M \quad (I\m
f)(x)$ and $Q\m (f,g)(x)$ expand to asymptotic series in
$1/m$ as
$m\to +\infty$. $\F\Bigl(\bigl(I\m f\bigr)(x)\Bigr)
=\bigl(If\bigr)(x)$ and $\F\bigl(Q\m
(f,g)(x)\bigr)=Q(f,g)(x)$,
where $I$ and $Q$ are the formal Berezin
transform and the formal twisted product
corresponding to the star-product with separation of
variables $\star$ on $(M,\omega_{-1})$ whose dual
star-product $\tilde\star$ on $(M,-\omega_{-1})$ is
parametrized by the formal form
$\tilde\omega=-(1/\nu)\omega_{-1}+\omega_{can}.$
\end{theorem}

\begin{remark}
As shown in \cite{Schlbia98} we have the following chain of
inequalities
\begin{equation}\label{E:symchain}
|I^{(m)}(f)|_{\infty}=
|\sigma(T_f^{(m)})|_\infty\quad\le\quad ||T_f^{(m)}||\quad
\le \quad|f|_\infty\ .
\end{equation}
Here $||..||$ denotes the operator norm with respect to the
norm of the sections of $L^m$ and $|..|_\infty$ the sup-norm
on
$C^{\infty}(M)$.
Choose as $x_e\in M$ a point with $|f(x_e)|=|f|_\infty$.
From \refT{central} and the fact that the formal Berezin
transform has
as leading term the identity   it follows that
$\ |(I^{(m)}f)(x_e)-f(x_e)|\le A/m\ $ with a suitable constant
$A$.
This implies
$\ \left| |f(x_e)|-|(I^{(m)}f)(x_e)| \right| \le A/m\ $
and hence
\begin{equation}\label{E:absch}
|f|_\infty-\frac Am=|f(x_e)|-\frac Am\quad\le\quad
|(I^{(m)}f)(x_e)|\quad\le\quad |(I^{(m)}f)|_\infty\ .
\end{equation}
Putting \refE{symchain} and \refE{absch} together we obtain
\begin{equation}\label{E:thma}
|f|_\infty-\frac Am\quad\le\quad ||T_f^{(m)}||\quad\le\quad
|f|_\infty\ .
\end{equation}
This provides another proof of \cite{BMS}, Theorem 4.1.
\end{remark}

\section{The identification of the Berezin-Toeplitz
star-product}\label{S:star}
In this section $\star$ will denote the
star-product with separation of variables on
$(M,\omega_{-1})$ whose dual $\tilde\star$ is
the star-product with separation of variables on
$(M,-\omega_{-1})$ parametrized by the formal
form $\tilde\omega=-(1/\nu)\omega_{-1}+
\omega_{can}$. 

Let $I=1+\nu I_1+\nu^2 I_2+\dots$ and $Q=Q_0+\nu
Q_1+\dots$ denote the formal Berezin transform
and the formal twisted product corresponding
to $\star$. \refT{central} asserts
that for given $f,g\in C^\infty(M),\ r\in{\bf
N},\ x\in M$ there exist constants $A,B$ such
that for sufficiently big values of $m$ the
following inequalities hold:

\begin{equation}\label{E:iii}
\left|\bigl(I^{(m)}f\bigr)(x)-\sum_{i=0}^{r-1} \frac{1}{m^i}
I_i(f)(x)\right|\leq \frac{A}{m^r},
\end{equation}
\begin{equation}\label{E:qqq}
\left|Q^{(m)}(f,g)(x)-\sum_{i=0}^{r-1} \frac{1}{m^i}
Q_i(f,g)(x)\right|\leq \frac{B}{m^r}.
\end{equation}

It was proved in \cite{Schlbia95},\cite{Schldef} that
Berezin-Toeplitz quantization on a compact
K\"ahler manifold $M$ gives rise to a
star-product on $M$. This star-product
$\star^{BT}$ is given by a sequence of bilinear
operators $\{C_k\},\ k\geq 0,$ on $C^\infty(M)$
satisfying the following conditions. For
$f,g\in C^\infty(M)$ and any $r\in{\bf N}$
there exists a constant $C$ such that
\begin{equation}\label{E:tt}
  \left\|T^{(m)}_fT^{(m)}_g
-T^{(m)}_{f\star_{[r]}g}\right\| \leq C/m^r,
\end{equation}
where $f\star_{[r]}g=\sum_{k=0}^{r-1} (1/m^k) 
C_k(f,g)$. The conditions \refE{tt} determine
the star-product $\star^{BT}$ uniquely.  We
call $\star^{BT}$ the Berezin-Toeplitz
star-product.

Recall that for $f,g\in C^\infty(M)\ 
\sigma(T^{(m)}_f)=I^{(m)}(f),\
\sigma(T^{(m)}_fT^{(m)}_g)=Q^{(m)}(f,g)$. 

Passing from operators to their covariant
symbols in \refE{tt} and using the inequality
$|\sigma(A)|\leq \|A\|$ we get that
\begin{equation}\label{E:qi}
  \left|Q^{(m)}(f,g)(x)-I^{(m)}(f\star_{[r]} g)(x)\right|\leq C/m^r.
\end{equation}

It follows from \refE{iii} that
\begin{equation}\label{E:ik}
\left|\frac{1}{m^k}I^{(m)}\bigl(C_k(f,g)\bigr)(x)-\sum_{i=0}^{r-k-1}
\frac{1}{m^{i+k}} I_i\bigl(C_k(f,g)\bigr)(x)\right|\leq
\frac{A_k}{m^r}.
\end{equation}
Summing up inequalities \refE{qqq} and \refE{ik} for
$k=0,1,\dots,r-1$, we obtain
that
\begin{multline}\label{E:long}
  \Big|\Bigl(Q^{(m)}(f,g)(x)-I^{(m)}\bigl(f\star_{[r]}g\bigr)(x)\Bigr)-\\
 \sum_{i=0}^{r-1}
\frac{1}{m^i}\Bigl(Q_i(f,g)(x)-\sum_{j+k=i}
I_j\bigl(C_k(f,g)\bigr)(x)\Bigr)\Big|\leq
\frac{D}{m^r}.
\end{multline}
for some constant $D$. It follows from
\refE{qi} and \refE{long}
that
$$
 \left|\sum_{i=0}^{r-1}
\frac{1}{m^i}\Bigl(Q_i(f,g)(x)-\sum_{j+k=i}
I_j\bigl(C_k(f,g)\bigr)(x)\Bigr)\right|\leq
\frac{E}{m^r},
$$
for some constant $E$, which infers that for
$i=0,1,\dots$
\begin{equation}\label{E:last}
 Q_i(f,g)=\sum_{j+k=i} I_j(C_k(f,g)).
\end{equation}

Equalities \refE{last} mean that
$Q(f,g)=I(f\star^{BT} g)$. Since $I$ is
invertible we immediately obtain that the
star-products $\star'$ and $\star^{BT}$
coincide. Thus the Berezin-Toeplitz deformation
quantization is completely identified as
the deformation quantization
with separation of variables on
$(\overline{M},\omega_{-1})$ whose
star-product $\star^{BT}$ is opposite to
$\tilde\star$. 

Using \refE{coh} we can calculate the
characteristic class $cl(\star^{BT})$
of the Berezin-Toeplitz
star-product $\star^{BT}$. 

It follows from \refE{coh} and
\refE{omega} that the characteristic class
of the star-product $\tilde\star$ equals
to $cl(\tilde\star)=
(1/i)\bigl(-[(1/\nu)\omega_{-1}]
+\varepsilon/2\bigr)$. It is easy to show
that the characteristic class of the
opposite star-product $\star'$ is equal to
$-cl(\tilde\star)$. Since
$\star^{BT}=\star'$, we finally get that
the characteristic class of the
Berezin-Toeplitz deformation quantization
is given by the formula $cl(\star^{BT})=
(1/i)\bigl([(1/\nu)\omega_{-1}]
-\varepsilon/2\bigr)$.

The characteristic class of the
Berezin-Toeplitz deformation
quantization was first calculated
by Eli Hawkins in \cite{Haw} by K-theoretic
methods.

As a concluding remark we would like to draw
the readers attention to the fact that the
classifying form $\omega$ of the star-product 
$\star$ is the formal object corresponding to the
asymptotic expansion as $m\to+\infty$
of the pullback $\omega^{(m)}$
of the Fubini-Study form on the projective space
$\P(H_m^*)$
via Kodaira embedding of $M$ into $\P(H_m^*)$. 
Here $H_m^*$ denotes the Hilbert
space dual to $H_m=\Gamma_{hol}(L^m)$ (see \refS{symbols}).
It was proved by
Zelditch \cite{Zel} that $\omega^{(m)}$ admits a complete
asymptotic expansion in $1/m$ as $m\to+\infty$. 
As an easy consequence of the results obtained
in this article one can show that 
$\F(\omega^{(m)})=\omega$.



\end{document}